\documentclass[11pt]{article}
\usepackage{amsmath}
\usepackage{amstext,amsbsy}
\usepackage{epsfig,multicol}
\usepackage{amsthm}
\usepackage{amssymb,latexsym}
\theoremstyle{plain}
\newtheorem{theorem}{Theorem}[section]
\newtheorem{lemma}[theorem]{Lemma}
\newtheorem{corollary}[theorem]{Corollary}
\newtheorem{proposition}[theorem]{Proposition}
\newtheorem{example}[theorem]{Example}

\theoremstyle{definition}
\newtheorem{definition}[theorem]{Definition}

\numberwithin{equation}{section}

\begin{document}

\title{Crossings and Nestings of Two Edges in
Set Partitions}
\author{Svetlana~Poznanovik$^{1}$,  Catherine~Yan$^{1,2,3}$ \\
\vspace{.3cm} \\
$^1$ Department of Mathematics, \\
Texas A\&M University, College Station, TX 77843. U. S. A. \\
\vspace{.3cm} \\
$^2$Center of Combinatorics, LPMC-TJKLC \\
Nankai University, Tianjin 300071, P. R. China \\
\vspace{.3cm}  \\
E-mail: $^1$spoznan@math.tamu.edu, $^2$cyan@math.tamu.edu
}
\date{}

\maketitle

\footnotetext[3]{The second author was supported in part by NSF grant
DMS-0245526 and DMS-0653846. }



\begin{abstract}

Let $\pi$ and $\lambda$ be two set partitions with the same number of
blocks. Assume $\pi$ is a partition of $[n]$.
 For any integer $l, m \geq 0$, let $\mathcal{T}(\pi, l)$ be
the set of partitions of $[n+l]$ whose restrictions to the last $n$
elements are isomorphic to $\pi$, and $\mathcal{T}(\pi, l, m)$ the
subset of $\mathcal{T}(\pi,l)$  consisting of those partitions
with exactly $m$ blocks. Similarly define  $\mathcal{T}(\lambda,
l)$ and  $\mathcal{T}(\lambda, l,m)$. We prove that if the
statistic $cr$ ($ne$), the number of crossings (nestings) of two
edges, coincides on the sets $\mathcal{T}(\pi, l)$ and
$\mathcal{T}(\lambda, l)$ for $l =0, 1$, then it coincides on
$\mathcal{T}(\pi, l,m)$ and   $\mathcal{T}(\lambda, l,m)$ for all
$l , m \geq 0$. These results extend the ones obtained by Klazar
on the distribution of crossings and nestings for matchings.
\end{abstract}

\section{Introduction and Statement of  Main Result}

In a recent paper \cite{klazar},  Klazar studied distributions of the
numbers of crossings and nestings of two edges in (perfect)
matchings. All matchings form an infinite tree $\mathcal{T}$ rooted at
the empty matching $\emptyset$, in which the children of a matching
$M$ are the matchings obtained from $M$ by adding to $M$ in all
possible ways a new first edge.
Given two matchings $M$ and $N$ on $[2n]$, Klazar decided
when the number of crossings (nestings) have identical distribution
on the levels of the two subtrees of $\mathcal{T}$ rooted at
$M$ and $N$.
%
%proved the following results for matchings \cite{klazar}:
%Let $M, N$ be two matchings on $[2n]$. For any integer $l \geq 0$,
%let $\mathcal{T}(M,l)$ be the set of matchings on $[2n+2l]$ which can
%be obtained  from $M$ by successively add $l$ times the first edge,
%and similarly for $mathcal{T}(N,l)$.
%Let  $s, t  \in \{cr, ne\}$, where $cr$ is the number of
%crossings in a matching and $ne$ is the number of nestings.  If the
%statistics $s$ and $t$ coincide, respectively, on the set of matchings
%$\mathcal{T}(M,l)$ and $\mathcal{T}(N,l)$ for $l=0, 1$, then they
%coincide on these sets for every $l \geq 0$.
In the last section
of \cite{klazar} Klazar raised the question as to apply the method to
other structures besides matchings.
In the present paper we consider set partitions,
which have a natural graphic representation by a set of arcs.
We establish the  Klazar-type results to the
distribution of crossings and nestings of two edges in set partitions.

Our approach follows  that of Klazar for matchings
\cite{klazar}, 
but is not a straightforward generalization. The structure of
set partitions is more complicated than that of matchings.
For example, partitions of $[n]$ may have different number of
blocks, while every matching of $[2n]$ has exactly $n$ blocks.
To get the results we  first defined an infinite tree
$\mathcal{T}(\Pi)$ on the set of
all set partitions, which, when restricted to matchings, is 
different than the one introduced  by Klazar.
We state our  main result in  Theorem 1.1, whose proof and
applications are given in Sections 2 and 3. 
Section 4 is devoted to the enumeration of the crossing/nesting similarity
classes. 
Though the ideas of the proofs are  similar to those in \cite{klazar}, 
in many places we have   to supply our own argument
to fit in the different structure,  and use a variety of combinatorial
structures, in particular, Motzkin paths, Charlier diagrams,
and binary sequences. 
We also analyze the joint generating function of the statistics 
$cr$ and $ne$ over partitions rooted at $\pi$, and derive a continued fraction 
expansion for general $\pi$.

We begin by introducing necessary notations. 
A (set) partition of $[n]=\{1, 2, \dots, n\}$ is a collection of disjoint nonempty subsets
of $[n]$, called blocks,  whose union is $[n]$.
A matching of $[2n]$ is a partition of $[2n]$ in $n$ two-element
blocks, which we also call {\em edges}.  If a partition $\pi$ has $k$ blocks,
we write $|\pi|=k$.
%      $\pi=B_1 - B_2 - \cdots \- B_k$, where
%$\min(B_1) < \mon(B_2) < \cdots < \min(B_k)$, and, within each block,
%the elements are listed in the numerical order.
A partition $\pi$ is
often represented as a graph on the vertex set $[n]$, drawn on a
horizontal line in the increasing order from left to right, whose
edge set consists of arcs connecting the elements of each block in
numerical order. We write an arc $e$ as a pair $(i,j)$ with $i<j$.

For  a partition $\pi$ of $[n]$, 
we say that the arcs $(i_1,j_1)$
and $(i_2,j_2)$ form a crossing if $i_1 < i_2 < j_1 < j_2$, and
they form a nesting if $i_1 <i_2 < j_2 <j_1$. By $cr(\pi)$ (resp.
$ne(\pi)$), we denote the number of crossings (resp. nestings) of
$\pi$. The distribution of the statistics $cr$ and $ne$ on matchings
has been studied in a number of articles, 
including \cite{CDDSY,klazar,Riordan,deSC,touchard}, to list a few.
The symmetry of $cr$ and $ne$ for set partitions was established 
in \cite{kz}. 
In this paper we investigate the distribution of the statistics 
$cr(\pi)$ and $ne(\pi)$ over the partitions of $[n]$ with a prefixed 
restriction to the last $k$ elements.

\begin{figure}[ht]
\begin{center}
\includegraphics[width=11cm]{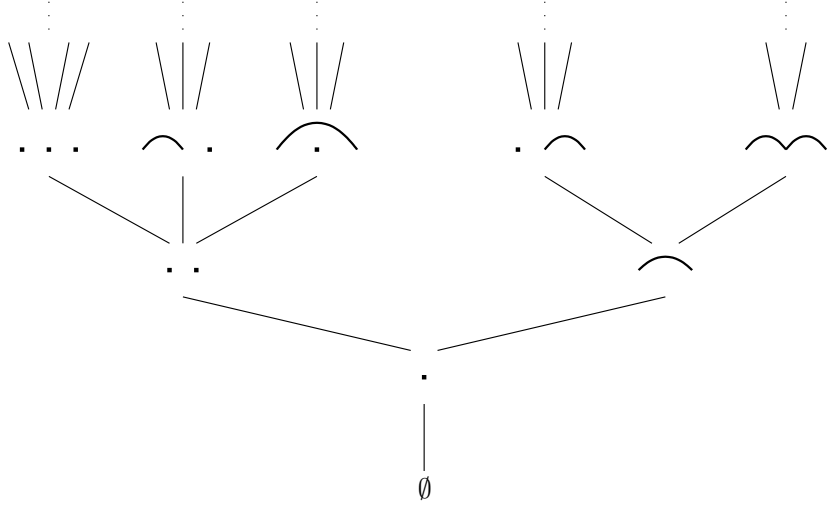}
\end{center}
\caption{The tree of partitions $\mathcal{T}(\Pi)$}
\label{tree_of_partitions}
\end{figure}

Denote by $\Pi_n$ the set of all partitions of $[n]$, 
and by $\Pi_{n,k}$ the set of partitions of $[n]$ with $k$ blocks. 
For $n=0$, $\Pi_0$ contains the empty partition. 
Let 
$\Pi=\cup_{n=0}^\infty \Pi_n= \cup_{n=0}^\infty \cup_{k \leq n} \Pi_{n,k}$. 
We define {\em the tree $\mathcal{T}(\Pi)$ of partitions} as a rooted tree
whose nodes are partitions such that:
\begin{enumerate}
\item The root is the empty partition; 
\item The partition $\pi $
of $[n+1]$ is a child of $\lambda$, a partition of [n], if and
only if the restriction of $\pi$ on $\{2,\dots,n+1 \}$ is
order-isomorphic to $\lambda$.
\end{enumerate}
See Figure \ref{tree_of_partitions} for an illustration of
$\mathcal{T}(\Pi)$. 

Observe that if $\lambda$ is a
partition of $[n]$ with $|\lambda|=k$, then $\lambda$ has $k+1$
children in $\mathcal{T}(\Pi)$. Let $B_{1},\dots,B_{k}$ be the blocks of $\lambda$
ordered in increasing order with respect to their minimal
elements. 
 For a set $S$, let $S+1=\{a+1: a \in S\}$. 
We denote the children of $\lambda$ by
$\lambda^{0},\lambda^{1},\dots,\lambda^{k}$ as follows:  
$\lambda^{0}$ is a partition of $[n+1]$ with $k+1$ blocks,
$$\lambda^{(0)}=\{ \{1\},B_{1}+1,\dots,B_{k}+1\};
$$
for $1 \leq i \leq k$, $\lambda^{i}$ is a partition of 
$[n+1]$ with $k$ blocks,
$$
 \lambda^{(i)}=\{ \{1\}
\cup (B_{i}+1),B_{1}+1,\dots,B_{i-1}+1,B_{i+1}+1,\dots,B_{k}+1 \}.
$$

For a partition $\lambda$, let $\mathcal{T}(\lambda)$ denote the
subtree of $\mathcal{T}(\Pi)$ rooted at $\lambda$, and let
$\mathcal{T}(\lambda,l)$ be the set of all partitions at the
$l$-th level of $\mathcal{T}(\lambda)$. $\mathcal{T}(\lambda,l,m)$
is the set of all partitions on the $l$-th level of
$\mathcal{T}(\lambda)$ with $m$ blocks. Note that
$\mathcal{T}(\lambda,l,m) \neq \emptyset$ if and only if $k \leq m
\leq k+l$.

Let $G$ be an abelian group and $\alpha$, $\beta$ two elements
in $G$.  Consider the statistics $s_{\alpha,\beta}: \Pi
\rightarrow G$ given by
$s_{\alpha,\beta}(\lambda)=cr(\lambda)\alpha + ne(\lambda)\beta$.
In \cite{klazar}, Klazar defines a tree of matchings and shows
that for two matchings $M$ and $N$, if the statistic
$s_{\alpha,\beta}$ coincides at the first two levels of
$\mathcal{T}(M)$ and $\mathcal{T}(N)$ then it coincides at all
levels, and similarly for the pair of statistics
$s_{\alpha,\beta},s_{\beta,\alpha}$. In this article we prove that
in the tree of partitions defined  above, the same results hold.
Precisely,

%%%%%%%%%%%%%%%%%%%%%%%%%%%%%%%%%%%%%%%%%%%%%%%%%%%%%%%%%%%%%%%%%%
% Main theorem                                                   %
%%%%%%%%%%%%%%%%%%%%%%%%%%%%%%%%%%%%%%%%%%%%%%%%%%%%%%%%%%%%%%%%%%
\begin{theorem}\label{T:main}
Let $\lambda, \pi \in \mathcal{T}(\Pi)$ be two non-empty partitions, and
 $s_{\alpha, \beta}(T) $ be the multiset containing
$\{ s_{\alpha, \beta}(t): t \in T\}$. We have
\begin{enumerate}
\item[(a)] If $s_{\alpha,
\beta}(\mathcal{T}(\lambda,l))=s_{\alpha,
\beta}(\mathcal{T}(\pi,l))$ for $l=0,1$ then \\
$s_{\alpha, \beta}(\mathcal{T}(\lambda,l,m))=s_{\alpha,
\beta}(\mathcal{T}(\pi,l,m))$ for all $l, m \geq 0$.

\item[(b)]If $s_{\alpha, \beta}(\mathcal{T}(\lambda,l))=s_{\beta,
\alpha}(\mathcal{T}(\pi,l))$ for $l=0,1$ then \\
$s_{\alpha, \beta}(\mathcal{T}(\lambda,l,m))=s_{\beta,
\alpha}(\mathcal{T}(\pi,l,m))$ for all $l,m \geq 0$.
\end{enumerate}
In other words, if the statistic $s_{\alpha,\beta}$ coincides on the
first two levels of the trees $\mathcal{T}(\lambda)$ and
$\mathcal{T}(\pi)$ then it coincides on $\mathcal{T}(\lambda,l,m)$
and $\mathcal{T}(\pi,l,m)$ on all levels, and similarly for the
pair of statistics $s_{\alpha,\beta},s_{\beta,\alpha}$.
\end{theorem}
Note that the conditions of Theorem 1.1 imply that $\lambda$ and
$\pi$ have the same number of blocks. But they are not necessarily
partitions of the same $[n]$. 

At the end of the introduction 
we would like to point out the major differences between the
structure of crossing and nesting  of set partitions and 
that of matchings. 
\begin{enumerate}
\item  The tree of partitions $\mathcal{T}(\Pi)$ and the tree of
  matchings are different. In $\mathcal{T}(\Pi)$, children of a
  partition $\pi$ is obtained by adding a new  vertex, instead of
 adding a first edge. Hence Klazar's  tree of matchings is not a
 sub-poset of $\mathcal{T}(\Pi)$. The definition of $\mathcal{T}(\Pi)$
  allows us to define the analogous  operators $R_{\alpha, \beta, i}$,
  as in \cite[\S 2]{klazar}.
  Since some descendants of $\pi$ are obtained by
adding isolated points, we need to introduce an extra  operator $M$,
(see Definition 2.2), and supply some new arguments to work with our
structure and $M$.   
\item  The type of a matching is encoded by a Dyck path, while for
set partitions, the corresponding structure is restricted bicolored
Motzkin paths(RBM), (c.f. Section 3).  
\item   The enumeration of crossing/nesting similarity classes is
different. A crossing-similarity class is determined
by a value $cr(M)$ ($cr(\pi)$) and a composition  $(a_1, a_2, \dots, a_m)$ of
$n$. 
For matchings  $cr(M)$ can be any integer between $0$ and 
$1+a_2+ 2a_3+\cdots + (m-1)a_m$. 
But for partitions the possible value of $cr(\pi)$ depends only on
$m$, but not $a_i$'s. 

In  matchings there is a bijection
between the set of nesting sequences of matchings of $[2n]$
 and the set of Dyck paths $\mathcal{D}(n)$. 
There is no analogous result   between set partitions and restricted bicolored 
Motzkin paths. 
\item  For matchings every nesting-similarity class is a subset of 
a crossing-similarity class. This is not true for set partitions. 
\end{enumerate}

%%%%%%%%%%%%%%%%%%%%%%%%%%%%%%%%%%%%%%%%%%%%%%%%%%%%%%%%%%%%%%%%%%%%

\section{Proof of Theorem~\ref{T:main}}\label{S:main}

%%%%%%%%%%%%%%%%%%%%%%%%%%%%%%%%%%%%%%%%%%%%%%%%%%%%%%%%%%%%%%%%%%%

Throughout this article we will generally adapt Klazar's notation on
multisets.
Formally a multiset is a pair $(A,m)$, where $A$ is a set, called the
underlying set, and $m:A \rightarrow \mathbb{N}$ is a mapping that
determines the multiplicities of the elements of $A$.
We often write multisets by repeating
the elements according to their multiplicities. 

For a map $f:X \rightarrow Y$ and $Z \subset X$,
let $f(Z)$ denote the multiset whose underlying set is
$\{f(z) :z \in Z\}$ and in which each element $y$ appears with
 multiplicity equal to the cardinality of
the set $\{z: z \in Z \text{  and  } f(z)=y\}$.
$\mathcal{S}(X)$  denotes the set of all finite multisets with
elements in the set $X$. Any
function $f:X \rightarrow S(Y)$ naturally extends to
$f:\mathcal{S}(X) \rightarrow \mathcal{S}(Y)$
by $f(Z)=\bigcup_{z \in Z}\{f(z)\}$, where
$\bigcup$ is union of multisets (the multiplicities of elements
are added).

For each $b_{i}=\min{B_{i}}$ of $\lambda$ define $u_{i}(\lambda)$
to be the number of edges $(p,q)$ such that $p< b_i <q$  and
$v_{i}(\lambda)$ to be the number of edges $(p,q)$ such that $p<q<
b_{i}$. They satisfy the obvious recursive relations
%%%%%%%%%%%%%%%%%%%%%%%%%%%%%%%%%%%%%%%%%%%%%%%%%%%%%%%%%%%%%%%%%%%
% Recursive relations for u and v                                 %
%%%%%%%%%%%%%%%%%%%%%%%%%%%%%%%%%%%%%%%%%%%%%%%%%%%%%%%%%%%%%%%%%%%
\begin{align}
u_{i}(\lambda^{0})&=
\begin{cases}
    0 &\text{if $i=1$} \\
    u_{i-1}(\lambda) &\text{if $2 \leq i \leq k+1$}
\end{cases} \label{E:uvrec1}
\\
v_{i}(\lambda^{0})&=
\begin{cases}
    0 &\text{if $i=1$} \\
    v_{i-1}(\lambda) &\text{if $2 \leq i \leq k+1$}
\end{cases}\label{E:uvrec2}
\\
u_{i}(\lambda^{j})&=
\begin{cases}
    0  &\text{if $i=1$} \\
    u_{i-1}(\lambda)+1 &\text{if $2 \leq i\leq j$}\\
    u_{i}(\lambda) &\text{if $j+1 \leq i \leq k$}
\end{cases} \label{E:uvrec3}
\\
v_{i}(\lambda^{j})&=
\begin{cases}
    0 &\text{if $i=1$} \\
    v_{i-1}(\lambda) &\text{if $2 \leq i\leq j$}\\
    v_{i}(\lambda)+1 &\text{if $j+1 \leq i \leq k$}
\end{cases}\label{E:uvrec4}
\end{align}
for $j=1, \dots , k$,  where $k=|\lambda| \geq 1 $.
%%%%%%%%%%%%%%%%%%%%%%%%%%%%%%%%%%%%%%%%%%%%%%%%%%%%%%%%%%%%%%%%%%%
For the statistics 
 $s_{\alpha,\beta}: \Pi \rightarrow G$ defined by
$s_{\alpha,\beta}(\lambda)=cr(\lambda)\alpha + ne(\lambda)\beta$,
we have that
\begin{align}
s_{\alpha,\beta}(\lambda^{0})&=s_{\alpha,\beta}(\lambda^{1})=s_{\alpha,\beta}(\lambda),
\label{E:first}
\\
s_{\alpha,\beta}(\lambda^{j})&=s_{\alpha,\beta}(\lambda)+u_{j}(\lambda)\alpha
+ v_{j}(\lambda)\beta,    \qquad j \geq 1. \label{E:second}
\end{align}
For simplicity, we will write $\lambda ^{ij}$ for $(\lambda^i)^j$.

%%%%%%%%%%%%%%%%%%%%%%%%%%%%%%%%%%%%%%%%%%%%%%%%%%%%%%%%%%%%%%%%%%%
% Lemma : Recurrence for s_alpha,beta                             %
%%%%%%%%%%%%%%%%%%%%%%%%%%%%%%%%%%%%%%%%%%%%%%%%%%%%%%%%%%%%%%%%%%%
\begin{lemma} \label{L:salbe}
For $|\lambda | \geq 1$,
\begin{align}
s_{\alpha,\beta}(\lambda^{0 j}))&=
\begin{cases}
s_{\alpha,\beta}(\lambda) &\text{if  $j=0,1$}\\
s_{\alpha,\beta}(\lambda^{j-1}) &\text{if  $j \geq 2$},\\
\end{cases}\label{E:salbe1}
\end{align} 
and for $i \geq 1$, 
\begin{align} 
 s_{\alpha,\beta}(\lambda^{i j})&=
\begin{cases}
s_{\alpha,\beta}(\lambda^{i}) &\text{if  $j=0,1$}\\
s_{\alpha,\beta}(\lambda^{i})+s_{\alpha,\beta}(\lambda^{j-1})-s_{\alpha,\beta}(\lambda^{1})+\alpha &\text{if $2 \leq j \leq i$}\\
s_{\alpha,\beta}(\lambda^{i})+s_{\alpha,\beta}(\lambda^{j})-s_{\alpha,\beta}(\lambda^{1})+\beta
&\text{if  $j \geq i+1$.}
\end{cases}\label{E:salbe2}
\end{align}
\end{lemma}

\begin{proof}
We first show \eqref{E:salbe2}. The first line in \eqref{E:salbe2}
follows directly from \eqref{E:first}. For the other two,
\begin{align*}
s_{\alpha,\beta}(\lambda^{i
j})&=s_{\alpha,\beta}(\lambda^{i})+u_{j}(\lambda^{i})\alpha +
v_{j}(\lambda^{i})\beta \\ \notag
 &=
\begin{cases}
s_{\alpha,\beta}(\lambda^{i})+u_{j-1}(\lambda)\alpha + \alpha +
v_{j-1}(\lambda)\beta &\text{if $2 \leq j \leq i$} \\
s_{\alpha,\beta}(\lambda^{i})+u_{j}(\lambda)\alpha +
v_{j}(\lambda)\beta + \beta &\text{if $j \geq i+1$}
\end{cases}\\ \notag
&=
 \begin{cases}
s_{\alpha,\beta}(\lambda^{i})+s_{\alpha,\beta}(\lambda^{j-1})-s_{\alpha,\beta}(\lambda^{1})+\alpha &\text{if $2 \leq j \leq i$}\\
s_{\alpha,\beta}(\lambda^{i})+s_{\alpha,\beta}(\lambda^{j})-s_{\alpha,\beta}(\lambda^{1})+\beta
&\text{if  $j \geq i+1$.}
\end{cases}\notag
\end{align*}
The first and third equality follow from \eqref{E:second} and the
second one follows from \eqref{E:uvrec3} and \eqref{E:uvrec4}.
Similarly, \eqref{E:salbe1} follows from \eqref{E:uvrec1},
\eqref{E:uvrec2}, \eqref{E:first}, and \eqref{E:second}.
\end{proof}
%%%%%%%%%%%%%%%%%%%%%%%%%%%%%%%%%%%%%%%%%%%%%%%%%%%%%%%%%%%%%%%%%%%%%%%%%%%%%%%%%%
To each partition $\lambda$ with $k$ blocks, ($k \geq 1$),  we associate a
sequence
\[
seq_{\alpha,\beta}(\lambda):=
s_{\alpha,\beta}(\lambda^{1})s_{\alpha,\beta}(\lambda^{2})\dots
s_{\alpha,\beta}(\lambda^{k})
\] 
The sequence
$seq_{\alpha,\beta}(\lambda)$ encodes the information about the
distribution of $s_{\alpha,\beta}$ on the children of $\lambda$ in
$\mathcal{T}(\Pi)$, in which 
$s_{\alpha,\beta}(\lambda^{1})$  plays a special role when 
we analyze the change of $seq_{\alpha,\beta}(\lambda)$ below .
This is due to the fact that 
$s_{\alpha,\beta}(\lambda^{1})$ carries information about $\lambda$
and  two children of $\lambda$, namely,  $\lambda ^0$ and $\lambda^1$.

For an abelian group $G$, let $G^*_l$  denote the set of finite
sequences of length $l$ over $G$, and $G^*=\cup_{l \geq 1} G^*_l$.
If $u=x_1x_2 \dots x_k \in G^*$ and $y \in G$,
then we use the convention that the sequence $(x_1+y)(x_2+y) \dots (x_k+y)$
is denoted  by $x_1x_2 \dots x_k +y$.
% Motivated by the results in Lemma~\ref{L:salbe}, we introduce the following

%%%%%%%%%%%%%%%%%%%%%%%%%%%%%%%%%%%%%%%%%%%%%%%%%%%%%%%%%%%%%%%%%%%%%%%%%%%%%%%
% Definition of R_alpha,beta                                                  %
%%%%%%%%%%%%%%%%%%%%%%%%%%%%%%%%%%%%%%%%%%%%%%%%%%%%%%%%%%%%%%%%%%%%%%%%%%%%%%%
\begin{definition}
For $\alpha, \beta \in G$ and $i \geq 1$, define  $R_{\alpha,
\beta,i}:G^{*}_l \rightarrow G^{*}_l$, ($ i \leq l$)  by setting
\[
R_{\alpha,\beta,i}(x_{1}x_{2}\dots x_{l})=x_{i}(x_{1} \dots x_{i-1} +
(x_{i}-x_{1}+\alpha))(x_{i+1} \dots x_{l}+(x_{i}-x_{1}+\beta))
\]
and $R_{\alpha, \beta}:G^{*} \rightarrow S(G^{*})$ by setting
\[
R_{\alpha, \beta}(x_{1}x_{2}\dots x_{l})=\{R_{\alpha,
\beta,i}(x_{1}x_{2}\dots x_{l}):1\leq i \leq l\}.
\]
In addition, define $M:G^{*}\rightarrow G^{*}$ by setting
\[
M(x_{1}x_{2}\dots
x_{l})=x_{1}x_{1}x_{2}\dots x_{l}.
\]
\end{definition}
%%%%%%%%%%%%%%%%%%%%%%%%%%%%%%%%%%%%%%%%%%%%%%%%%%%%%%%%%%%%%%%%%%%%%%%%%%%%%%
Lemma~\ref{L:salbe} immediately implies that
\begin{align*}
seq_{\alpha,\beta}(\lambda^{0})&= M(seq_{\alpha,\beta}(\lambda)),
\\
seq_{\alpha,\beta}(\lambda^{i})&=
R_{\alpha,\beta,i}(seq_{\alpha,\beta}(\lambda)),
\text{\qquad for } 1 \leq i \leq |\lambda|.
\end{align*}
For $l \geq 0$, let
$E_{\alpha,\beta}(\lambda,l,m)=\{seq_{\alpha,\beta}(\mu): \mu \in
\mathcal{T}(\lambda,l,m)\}$, the multiset of sequences
$seq_{\alpha,\beta}(\mu)$ associated to partitions $\mu \in
\mathcal{T}(\lambda,l,m)$. Then for $l \geq 1$,
\begin{equation}\label{E:levelseq}
E_{\alpha,\beta}(\lambda,l,m)=
R_{\alpha,\beta}(E_{\alpha,\beta}(\lambda,l-1,m)) \cup
M(E_{\alpha,\beta}(\lambda,l-1,m-1)).
\end{equation}
Next, we define an auxiliary function $f$ which reflects the
change of the statistic $s_{\alpha,\beta}$ along $\mathcal{T}(\Pi)$.
Then we prove two general properties of $f$ and
use these properties to prove
Theorem~\ref{T:main}.
%%%%%%%%%%%%%%%%%%%%%%%%%%%%%%%%%%%%%%%%%%%%%%%%%%%%%%%%%%%%%%%%%%%
% Definition of f                                                 %
%%%%%%%%%%%%%%%%%%%%%%%%%%%%%%%%%%%%%%%%%%%%%%%%%%%%%%%%%%%%%%%%%%%
For an integer $r \geq 0$ and $\gamma \in G$, the function
$f:G^*\rightarrow S(G)$ is defined by
\[
f_{\gamma}^{r}(x_{1}x_{2}\dots x_{l}):=\{x_{a_1}+x_{a_2}+ \cdots +
x_{a_r}-(r-1)x_{1}+\gamma: 1<a_1<a_2<\cdots < a_r \leq l\} \] In
particular,
\begin{align*}
f_{0}^{0}(x_{1}x_{2}\dots x_{l})&=\{x_{1}\},\\
f_{0}^{1}(x_{1}x_{2}\dots x_{l})&=\{x_{2},\dots ,x_{l}\}.
\end{align*}
%%%%%%%%%%%%%%%%%%%%%%%%%%%%%%%%%%%%%%%%%%%%%%%%%%%%%%%%%%%%%%%%%%%%
% First lemma for f                                                %
%%%%%%%%%%%%%%%%%%%%%%%%%%%%%%%%%%%%%%%%%%%%%%%%%%%%%%%%%%%%%%%%%%%%
\begin{lemma}\label{L:big}

Let $X,Y \in S(G^*)$ be two multisets such that
$f_{\gamma}^{r}(X)=f_{\gamma}^{r}(Y)$ for every $r \geq 0$ and
$\gamma \in G$. Then

\begin{itemize}
\item[(a)] $f_{\gamma}^{r}(M(X))=f_{\gamma}^{r}(M(Y))$,

\item[(b)] $f_{\gamma}^{r}(R_{\alpha,
\beta}(X))=f_{\gamma}^{r}(R_{\alpha, \beta}(Y))$,

\item[(c)] $f_{\gamma}^{r}(R_{\alpha,
\beta}(X))=f_{\gamma}^{r}(R_{\beta, \alpha}(Y))$,
\end{itemize}
for every $r \geq 0$ and $\gamma \in G$.
\end{lemma}

\begin{proof}
(a)  
The elements in $f_{\gamma}^{r}(M(X))$ have the form
$y_{a_1}+y_{a_2}+ \cdots + y_{a_r}-(r-1)y_{1}+\gamma$ for some
$y_{1}y_{2}\dots y_{l+1} \in M(X)$, where $y_{1}y_{2}\dots
y_{l+1}=x_{1}x_{1}x_{2}\dots x_{l}$ for some $x_{1}x_{2}\dots
x_{l}\in X$. For $r=0$,
\[
f_{\gamma}^{0}(M(X))=\{y_1 + \gamma:y_{1}y_{2}\dots y_{l+1} \in
M(X)\}=\{x_{1}+\gamma: x_{1}x_{2}\dots x_{l}\in X\}=
f_{\gamma}^{0}(X).
\] 
Hence
$f_{\gamma}^{0}(X)=f_{\gamma}^{0}(Y)$ implies
$f_{\gamma}^{0}(M(X))=f_{\gamma}^{0}(M(Y))$.

For $r \geq 1$,
divide the multiset
$f_{\gamma}^{r}(M(X))$ into two disjoint multisets, 
\[
A= \{y_{a_1}+y_{a_2}+ \cdots + y_{a_r}-(r-1)y_{1}+\gamma:
y_{1}y_{2}\dots y_{l+1} \in M(X),\; a_1=2\}
\]
and
\[
B=\{y_{a_1}+y_{a_2}+ \cdots + y_{a_r}-(r-1)y_{1}+\gamma:
y_{1}y_{2}\dots y_{l+1} \in M(X),\; a_1>2\}.
\]
The
elements of $A$  can be written as
\begin{align*}
 y_{a_1}+y_{a_2}+ \cdots + y_{a_r}-(r-1)y_{1}+\gamma &=
x_{1}+y_{a_2}+ \cdots + y_{a_r}-(r-1)x_{1}+\gamma \\
&= y_{a_2}+ \cdots + y_{a_r}-(r-2)x_{1}+\gamma \\
& =x_{a_2-1}+ \cdots + x_{a_r-1}-(r-2)x_{1}+\gamma.
\end{align*}
Since $a_2-1>a_1-1=1$, the multiset $A$ is equal to
$f_{\gamma}^{r-1}(X)$. The elements in $B$
can be written as 
\[ 
y_{a_1}+y_{a_2}+ \cdots +y_{a_r}-(r-1)y_{1}+\gamma
=x_{a_1-1}+x_{a_2-1}+ \cdots +
x_{a_r-1}-(r-1)x_{1}+\gamma. \] 
Since $a_1 \geq 3$, the indices on
the right-hand side run through all the increasing $r$-tuples
$1<a_1-1<a_2-1 < \cdots < a_r-1 \leq l$. Therefore, $B$
is equal to $f_{\gamma}^{r}(X)$. So,
\begin{equation} \label{E:mrec}
f_{\gamma}^{r}(M(X))=f_{\gamma}^{r-1}(X)\cup f_{\gamma}^{r}(X).
\end{equation}
By assumption  we have
\[ f_{\gamma}^{r}(M(X))=f_{\gamma}^{r-1}(X)\cup
f_{\gamma}^{r}(X)= f_{\gamma}^{r-1}(Y) \cup
f_{\gamma}^{r}(Y)=f_{\gamma}^{r}(M(Y)). 
\]
(c) We will prove only (c) because the proof of (b) is similar and
easier. Since $f_{\gamma}^{r}(X)$ is a translation of
$f_{0}^{r}(X)$ by $\gamma$, it is enough to prove the result for
$\gamma =0$ only. The elements of $f_{0}^{r}(R_{\alpha,
\beta}(X))$ have the form $y_{a_1}+y_{a_2}+ \cdots +
y_{a_r}-(r-1)y_{1}$, where $y_{1}y_{2}\dots y_{l} \in R_{\alpha,
\beta}(X)$ is equal to $x_{i}(x_{1} \dots x_{i-1} +
x_{i}-x_{1}+\alpha)(x_{i+1} \dots x_{l}+x_{i}-x_{1}+\beta)$ for
some $x_{1}x_{2}\dots x_{l}\in X$ and  $i \in [l]$.

For $0 \leq t \leq r$,  let
\begin{multline*}
C_{t,\alpha,\beta}(X)=\{y_{a_1}+y_{a_2}+ \cdots +
y_{a_r}-(r-1)y_{1}: \\
y_{1}y_{2}\dots y_{l} \in R_{\alpha,
\beta,i}(X) \; \mbox{and} \;a_{t} \leq i < a_{t+1}, \; \mbox{for
some} \; i \in [l]\}.
\end{multline*}
An element $y_{a_1}+y_{a_2}+ \cdots +
y_{a_r}-(r-1)y_{1}  \in C_{t,\alpha, \beta}(X)$ is equal to
\begin{eqnarray}
& & x_{a_1-1}+ \cdots x_{a_t-1} +t(x_{i}-x_{1}+\alpha) +
x_{a_{t+1}}+\cdots + x_{a_{r}}+(r-t)(x_{i}-x_{1}+\beta)
-(r-1)x_{i}  \notag \\ 
 &= & x_{a_1-1}+ \cdots
x_{a_t-1}+x_{i}+x_{a_{t+1}}+\cdots +
x_{a_{r}}-rx_{1}+t\alpha+(r-t)\beta. \label{E:type}
\end{eqnarray}
Again, we consider two cases, according to the value of $a_1$.
By~\eqref{E:type}, the submultiset
of $C_{t,\alpha, \beta}(X)$ for $a_1>2$ is equal to
$f^{r+1}_{t\alpha+(r-t)\beta}(X)$, and for $a_1=2$ the
corresponding submultiset is equal to
$f^{r}_{t\alpha+(r-t)\beta}(X)$. Therefore,
\begin{eqnarray}
\label{E:CX}
C_{t,\alpha, \beta}(X)=
f^{r+1}_{t\alpha+(r-t)\beta}(X) \cup
f^{r}_{t\alpha+(r-t)\beta}(X).
\end{eqnarray}
Similarly,
\begin{eqnarray}
\label{E:CY}
C_{t,\beta, \alpha}(Y)=
f^{r+1}_{t\beta+(r-t)\alpha}(Y) \cup
f^{r}_{t\beta+(r-t)\alpha}(Y).
\end{eqnarray} So,
\begin{align*}
f^{r}_{0}(R_{\alpha,\beta}(X))&= \bigcup_{t=0}^{r}C_{t,\alpha,
\beta}(X) \\
&=\bigcup_{t=0}^{r}f^{r+1}_{t\alpha+(r-t)\beta}(X)\cup
\bigcup_{t=0}^{r}f^{r}_{t\alpha+(r-t)\beta}(X)\\
&=\bigcup_{t=0}^{r}f^{r+1}_{t\alpha+(r-t)\beta}(Y)\cup
\bigcup_{t=0}^{r}f^{r}_{t\alpha+(r-t)\beta}(Y)\\
&=\bigcup_{t=0}^{r}f^{r+1}_{(r-t)\alpha+t\beta}(Y)\cup
\bigcup_{t=0}^{r}f^{r}_{(r-t)\alpha+t\beta}(Y)\\
&=\bigcup_{t=0}^{r}C_{t,\beta,
\alpha}(Y) \\
&=f^{r}_{0}(R_{\beta,\alpha}(Y)).
\end{align*}
The second and fifth equality follow from \eqref{E:CX} and
\eqref{E:CY} respectively.  The third equality follows from the
assumption of the lemma, while the fourth equality is just a
reordering of the unions.
\end{proof}
%%%%%%%%%%%%%%%%%%%%%%%%%%%%%%%%%%%%%%%%%%%%%%%%%%%%%%%%%%%%%%%%%%%%%%%%%%%
% The small lemma for f                                                   %
%%%%%%%%%%%%%%%%%%%%%%%%%%%%%%%%%%%%%%%%%%%%%%%%%%%%%%%%%%%%%%%%%%%%%%%%%%%
\begin{lemma}\label{L:small}

If $X,Y \in S(G^*)$  are one-element sets such that
$f_{0}^{0}(X)=f_{0}^{0}(Y)$ and $f_{0}^{1}(X)=f_{0}^{1}(Y)$, then
$f_{\gamma}^{r}(X)=f_{\gamma}^{r}(Y)$ for every $r \geq 0$ and
$\gamma \in G$.
\end{lemma}

\begin{proof} We need to prove that if $u,v \in G^{*}$ are two
sequences beginning with the same term and having equal numbers of
occurrences of each $g \in G$, then
$f_{\gamma}^{r}(u)=f_{\gamma}^{r}(v)$ for every $r \geq 0$ and
$\gamma \in G$. It suffices to prove the statement for $\gamma=0$,
because $f_{\gamma}^{r}(u)$ is a translation of $f_{0}^{r}(u)$ by
$\gamma$. Let $u=u_{1}\dots u_{l}$ and $v=v_{1}\dots v_{l}$. Since
$u_{1}=v_{1}$, it suffices to prove that the multisets
$\{u_{a_{1}}+u_{a_{2}}+ \cdots +u_{a_{r}}: 1 < a_{1} < a_{2} <
\dots < a_{r} \leq l \}$ and $\{v_{a_{1}}+v_{a_{2}}+ \cdots
+v_{a_{r}}: 1 < a_{1} < a_{2} < \dots < a_{r} \leq l \}$ are
equal. That is clear because $\{u_{2},\dots,u_{l}\}$ and
$\{v_{2},\dots,v_{l}\}$ are equal as multisets.
\end{proof}
%%%%%%%%%%%%%%%%%%%%%%%%%%%%%%%%%%%%%%%%%%%%%%%%%%%%%%%%%%%%%%%%%%%%%%%%%%%%%
% Proof of main theorem                                                     %
%%%%%%%%%%%%%%%%%%%%%%%%%%%%%%%%%%%%%%%%%%%%%%%%%%%%%%%%%%%%%%%%%%%%%%%%%%%%%
\begin{proof}[\textbf{Proof of Theorem~\ref{T:main} }]

We prove (b), the proof of (a) is similar. First, we  prove by
induction on $l$ that
\begin{eqnarray}\label{E:ind}
f_{\gamma}^{r}(E_{\alpha,\beta}(\lambda,l,m))=
f_{\gamma}^{r}(E_{\beta,\alpha}(\pi,l,m)) \text{\; for every  $r
\geq 0$ and $\gamma \in G$.}
 \end{eqnarray} 
Before we proceed with
the induction, it is useful to observe that the assumption \[
s_{\alpha, \beta}(\mathcal{T}(\lambda,l))=s_{\beta,
\alpha}(\mathcal{T}(\pi,l)) \text{\; for } l=0,1\] of
Theorem~(\ref{T:main}) (b) is equivalent to
\begin{eqnarray} \label{E:equiv}
s_{\alpha,\beta}(\mathcal{T}(\lambda,l,m))=s_{\beta,
\alpha}(\mathcal{T}(\pi,l,m)) \text{\; for $l=0,1$ and $k \leq m
\leq k+l$},
\end{eqnarray}
where $k = |\lambda|$. 
One direction is clear, the other one
follows from the following equations. 
\begin{align*}
s_{\alpha, \beta}(\mathcal{T}(\lambda,1,k+1))&=s_{\alpha,
\beta}(\mathcal{T}(\lambda,0,k))=s_{\alpha,
\beta}(\mathcal{T}(\lambda,0)),\\ 
s_{\alpha,\beta}(\mathcal{T}(\lambda,1,k))&=s_{\alpha,
\beta}(\mathcal{T}(\lambda,1)) \backslash s_{\alpha,
\beta}(\mathcal{T}(\lambda,0)),
\end{align*}
where $\backslash$ is the difference of multisets. For the same
reason the assumption of part (a) is equivalent to
\[s_{\alpha,\beta}(\mathcal{T}(\lambda,l,m))=s_{\alpha,
\beta}(\mathcal{T}(\pi,l,m)) \text{\; for $l=0,1$ and $k \leq m
\leq k+l$.}\]

Now we show \eqref{E:ind}. For $l=0$ we need to show
$f_{\gamma}^{r}(E_{\alpha,\beta}(\lambda,0,k))=
f_{\gamma}^{r}(E_{\beta,\alpha}(\pi,0,k))$. By Lemma~\ref{L:small}
we only need to check that $f_{0}^{0}(X)=f_{0}^{0}(Y)$ and
$f_{0}^{1}(X)=f_{0}^{1}(Y)$ for $X=\{seq_{\alpha,
\beta}(\lambda)\}$ and $Y=\{seq_{\beta,\alpha}(\pi)\}$. This
follows from \eqref{E:equiv}, because $f_{0}^{0}(X)=s_{\alpha,
\beta}(\mathcal{T}(\lambda,0,k))$ and $f_{0}^{1}(X)=s_{\alpha,
\beta}(\mathcal{T}(\lambda,1,k))$.

Suppose $f_{\gamma}^{r}(E_{\alpha,\beta}(\lambda,s,m))=
f_{\gamma}^{r}(E_{\beta,\alpha}(\pi,s,m))$ for all $0 \leq s <l$
and all $m$. Then using  ~\eqref{E:levelseq}, the induction
hypothesis, and Lemma~\ref{L:big} we have
\begin{align*}
f_{\gamma}^{r}(E_{\alpha,\beta}(\lambda,l,m))&=f_{\gamma}^{r}(R_{\alpha,\beta}(E_{\alpha,\beta}(\lambda,l-1,m)))
\cup f_{\gamma}^{r}(M(E_{\alpha,\beta}(\lambda,l-1,m-1)))\\
&=f_{\gamma}^{r}(R_{\beta,\alpha}(E_{\beta,\alpha}(\pi,l-1,m)))
\cup f_{\gamma}^{r}(M(E_{\beta,\alpha}(\pi,l-1,m-1)))\\
&=f_{\gamma}^{r}(E_{\beta,\alpha}(\pi,l,m)),
\end{align*}
and the induction is completed. 
This proves \eqref{E:ind}. 
Now
\[
s_{\alpha,\beta}(\mathcal{T}(\lambda,l,m))=f_{0}^{0}(E_{\alpha,\beta}(\lambda,l,m))
= f_{0}^{0}(E_{\beta,\alpha}(\pi,l,m))
=s_{\beta, \alpha}(\mathcal{T}(\pi,l,m)).
\]
\end{proof}

%%%%%%%%%%%%%%%%%%%%%%%%%%%%%%%%%%%%%%%%%%%%%%%%%%%%%%%%%%%%%%%%%%%%%%%%%%%%%%%%%%%%%%%
%%%%%%%%%%%%%%%%%%%%%%%%%%%%%%%%%%%%%%%%%%%%%%%%%%%%%%%%%%%%%%%%%%%%%%%%%%%%%%%%%%%%%%%
%                    Section: Some Applications and Examples                          %
%%%%%%%%%%%%%%%%%%%%%%%%%%%%%%%%%%%%%%%%%%%%%%%%%%%%%%%%%%%%%%%%%%%%%%%%%%%%%%%%%%%%%%%
%%%%%%%%%%%%%%%%%%%%%%%%%%%%%%%%%%%%%%%%%%%%%%%%%%%%%%%%%%%%%%%%%%%%%%%%%%%%%%%%%%%%%%%
\section{Applications and Examples}

As a direct corollary, we obtain a result of Kasraoui and Zeng
\cite[Eq.(1.6)]{kz}.

\begin{corollary}
The joint distribution of crossings and nestings of partitions is
symmetric i.e.
%\[|\{\lambda \in \Pi _n : cr(\lambda)=k,
%ne(\lambda)=m\}|=|\{\lambda \in \Pi _n : cr(\lambda)=m,
%ne(\lambda)=k\}|\]
\[
\sum_{\pi \in \Pi _n} p^{cr(\pi)}q^{ne(\pi)}=\sum_{\pi \in \Pi
_n} p^{ne(\pi)}q^{cr(\pi)}\]
\end{corollary}

\begin{proof}
Let $G=(\mathbb{Z} \oplus \mathbb{Z},+), \alpha =(1,0)$ and $\beta
=(0,1)$. The result follows from the second part of
Theorem~\ref{T:main} for $\lambda = \pi =\{\{1\}\}$.
\end{proof}
For a partition $\lambda$ we say that two edges form an alignment
if they neither form a crossing nor a nesting. The total number of
alignments in $\lambda$ is denoted by $al(\lambda)$. A stronger
result of Kasraoui and Zeng \cite[Eq.~(1.4)]{kz} can also be derived from
 Theorem~\ref{T:main}. 
\begin{corollary}
\[\sum_{\pi \in \Pi _n}{p^{cr(\pi)}q^{ne(\pi)}t^{al(\pi)}}=\sum_{\pi \in \Pi
_n}{p^{ne(\pi)}q^{cr(\pi)}t^{al(\pi)}}\]
\end{corollary}

\begin{proof}
Again we use $G=(\mathbb{Z} \oplus \mathbb{Z},+)$, $\alpha
=(1,0)$, $\beta =(0,1)$, and $\lambda = \pi =\{\{1\}\}$. Any
partition $\mu \in \Pi _n$ with $k$ blocks has $n-k$ edges. Hence
$cr(\mu)+ne(\mu)+al(\mu)= \binom{n-k}{2}$. The result follows from
the second part of Theorem~\ref{T:main}.
\end{proof}

\begin{corollary}
Let $\lambda$ and $\pi$ be two partitions of $[n]$ with same
number of blocks $k$. If the statistic $al$ is equidistributed on
the first two levels of $\mathcal{T}(\lambda)$ and
$\mathcal{T}(\pi)$, it is equidistributed on
$\mathcal{T}(\lambda,l,m)$ and $\mathcal{T}(\pi,l,m)$ for all $l,m
\geq 0$.
\end{corollary}

\begin{proof}
Again we use the identity $cr(\mu)+ne(\mu)+al(\mu)=
\binom{n-k}{2}$, which holds for any partition $\mu \in \Pi _n$
with $k$ blocks. Moreover, $al(\lambda)=al(\lambda ^{0})$.
Therefore the condition that the statistic $al$ is equidistributed
on the first two levels of $\mathcal{T}(\lambda)$ and
$\mathcal{T}(\pi)$ implies that the statistic $cr+ne$ is
equidistributed on $\mathcal{T}(\lambda,l,m)$ and
$\mathcal{T}(\pi,l,m)$ for all $l=0,1$ and all $m$. In other
words, if we set $G=\mathbb{Z}$ and $\alpha = \beta =1$ then the
the assumption of Theorem~\ref{T:main} is satisfied, and hence
$cr+ne$ is equidistributed on $\mathcal{T}(\lambda,l,m)$ and
$\mathcal{T}(\pi,l,m)$ for all $l,m \geq 0$. This, in return,
implies that $al$ is equidistributed on $\mathcal{T}(\lambda,l,m)$
and $\mathcal{T}(\pi,l,m)$ for all $l,m \geq 0$.
\end{proof}

\begin{example}
Let $\lambda=\{\{1,2,5\},\{3,4\}\}$ and
$\pi=\{\{1,2,4\},\{3,5\}\}$. There are as many partitions on $[n]$
with $m$ crossings and $l$ nestings which restricted to the last
five points form a partition isomorphic to $\lambda$ as there are
partitions of $[n]$ with $l$ crossings and $m$ nestings which
restricted to the last five points form a partition isomorphic to
$\pi$.
\end{example}
\begin{proof}
Set $G=(\mathbb{Z} \oplus \mathbb{Z},+), \alpha =(1,0)$ and $\beta
=(0,1)$, $s_{\alpha,\beta}=(cr,ne)$. The claim follows from part
(b) of Theorem~\ref{T:main} since
$s_{\alpha,\beta}(\lambda)=(0,1)=s_{\beta,\alpha}(\pi)$ and
$s_{\alpha,\beta}(\mathcal{T}(\lambda,1))=\{(0,1),(0,1),(1,2)\}=s_{\beta,\alpha}(\mathcal{T}(\pi,1))$
\end{proof}

\begin{example}
Let $\lambda=\{\{1,7\},\{2,6\},\{3,4\},\{5,8\}\}$ and
$\pi=\{\{1,8\},\{2,4\},\{3,6\},\{5,7\}\}$. There are as many
partitions on $[n]$ with $m$ crossings and $l$ nestings which
restricted to the last eight points form a partition isomorphic to
$\lambda$ as there are ones which restricted to the last eight
points form a partition isomorphic to $\pi$.

\end{example}
\begin{proof}
Again set $G=(\mathbb{Z} \oplus \mathbb{Z},+), \alpha =(1,0)$ and
$\beta =(0,1)$. Then $s_{\alpha,\beta}=(cr,ne)$. The claim follows
from part (a) of Theorem~\ref{T:main} since
\[
s_{\alpha,\beta}(\lambda)=(2,3)=s_{\alpha,\beta}(\pi)\]
 and
\[
s_{\alpha,\beta}(\mathcal{T}(\lambda,1))=\{(2,3),(2,3),(3,3),(4,3),(4,4)\}=s_{\alpha,\beta}(\mathcal{T}(\pi,1)).
\]
\end{proof}
%%%%%%%%%%%%%%%%%%%%%%%%%%%%%%%%%%%%%%%%%%%%%%%%%%%%%%%%%%%%%%%%%%%%%%%%%%%%%%%%%%%%%%%%%%%%%%%%%%%
%%%%%%%%%%%%%%%%%%%%%%%%%%%%%%%%%%%%%%%%%%%%%%%%%%%%%%%%%%%%%%%%%%%%%%%%%%%%%%%%%%%%%%%%%%%%%%%%%%%
%                                   Section: Crossing Classes                                     %
%%%%%%%%%%%%%%%%%%%%%%%%%%%%%%%%%%%%%%%%%%%%%%%%%%%%%%%%%%%%%%%%%%%%%%%%%%%%%%%%%%%%%%%%%%%%%%%%%%%
%%%%%%%%%%%%%%%%%%%%%%%%%%%%%%%%%%%%%%%%%%%%%%%%%%%%%%%%%%%%%%%%%%%%%%%%%%%%%%%%%%%%%%%%%%%%%%%%%%%

\section{Number of Crossing and Nesting-similarity Classes}

 In this section we consider equivalence relations $\sim_{cr}$ and
$\sim_{ne}$ on set partitions in the same way Klazar
defines them on matchings \cite{klazar} . We determine the number of
crossing-similarity classes in $\Pi_{n,k}$. 
For $\sim_{ne}$, we find a recurrence relation for the number
of nesting-similarity classes in $\Pi_{n,k}$, 
and compute the total number of such
classes in  $\Pi_n$. 

Define an equivalence relation $\sim_{cr}$ on $\Pi_n$: $\lambda
\sim_{cr} \pi$ if and only if
$cr(\mathcal{T}(\lambda,l,m))=cr(\mathcal{T}(\pi,l,m))$ for all
$l,m \geq 0$. 
The relation $\sim_{cr}$ partitions
$\Pi_{n,k}$ into equivalence classes.  Theorem~\ref{T:main}
implies that $\lambda \sim_{cr} \pi$ if and only if
$cr(\lambda)=cr(\pi)$ and
$f_0^1(seq_{1,0}(\lambda))=f_0^1(seq_{1,0}(\pi))$. Define
$crseq(\lambda)=seq_{1,0}(\lambda)-cr(\lambda)$. For the upcoming
computations it is useful to observe that $\lambda \sim_{cr} \pi$
if and only if $cr(\lambda)=cr(\pi)$ and
$f_0^1(crseq(\lambda))=f_0^1(crseq(\pi))$, i.e., $\lambda$ and
$\pi$ are equivalent if and only if they have the same number of
crossings and their sequences $crseq(\lambda)$ and $crseq(\pi)$
are equal as multisets. Denote the multiset consisting of the
elements of $crseq(\lambda)$ by $crset(\lambda)$.

Similarly, define $\lambda \sim_{ne} \pi$ if and only if
$ne(\mathcal{T}(\lambda,l,m))=ne(\mathcal{T}(\pi,l,m))$ for all
$l,m \geq 0$. Again, from Theorem~\ref{T:main} we have that
$\lambda \sim_{ne} \pi$ if and only if $ne(\lambda)=ne(\pi)$ and
$f_0^1(seq_{0,1}(\lambda))=f_0^1(seq_{0,1}(\pi))$. Since the
sequence $seq_{0,1}(\lambda)$ is nondecreasing, $\lambda \sim_{ne}
\pi$ if and only if $ne(\lambda)=ne(\pi)$ and
$seq_{0,1}(\lambda)-ne(\lambda)=seq_{0,1}(\pi)-ne(\pi)$. With the
notation at the beginning of Section~\ref{S:main},
$seq_{0,1}(\lambda)-ne(\lambda)=v_1 \dots v_k$. Denote this
sequence by $neseq(\lambda)$.

%%%%%%

A Motzkin path $M=(s_1,\dots,s_n)$ is a path from $(0,0)$ to
$(n,0)$ consisting of steps $s_i \in \{(1,1),(1,0),(1,-1)\}$ which
does not go below the $x$-axis. We say that the step $s_i$ is of
height $l$ if its left endpoint is at the line $y=l$. A restricted
bicolored Motzkin path is a Motzkin path with each horizontal step
colored red or blue which does not have a blue horizontal step
 of height 0. We will denote the steps $(1,1)$, $(1,-1)$, red $(1,0)$, and blue $(1,0)$ by
 NE (northeast), SE (southeast), RE (red east), and BE (blue east) respectively. The
set of all restricted bicolored Motzkin paths of length $n$ is
denoted by $RBM_n$. A Charlier diagram of length $n$ is a pair
$h=(M, \xi)$ where $M=(s_1,\dots,s_n) \in RBM_n$ and
$\xi=(\xi_1,\dots,\xi_n)$ is a sequence of integers such that
$\xi_i=1$ if $s_i$ is a NE or RE step, and $1 \leq \xi_i \leq l$
if $s_i$ is a SE or BE step of height $l$. $\Gamma_n$ will denote
the set of Charlier diagrams of length $n$.

It is well known that partitions are in one-to-one correspondence
with Charlier diagrams. Here we use two  maps
described in~\cite{kz}, which are based on similar constructions 
in  ~\cite{flajolet,viennot}.  
For our purpose, we reformulate the maps $\Phi_r,\Phi_l:\Gamma_n
\rightarrow \Pi_n$ as follows. Given $(M, \xi) \in \Gamma_n$,
construct $\lambda \in \Pi_n$ step by step. The path
$M=(s_1,\dots,s_n)$ determines the type of $\lambda$: $i \in [n]$
is

\begin{itemize}
\item[-] a minimal but not a maximal element of a block of
$\lambda$ (opener) if and only if $s_i$ is a NE step;

\item[-] a maximal but not a minimal element of a block of
$\lambda$ (closer) if and only if  $s_i$ is a SE step;

\item[-] both a minimal and a maximal element of a bock of
$\lambda$ (singleton) if and only if $s_i$ is a RE step;

\item[-] neither a minimal nor a maximal element of a block of
$\lambda$ (transient) if and only if  $s_i$ is a BE step.
\end{itemize}

To draw the edges in $\Phi_r((M,\xi))$, we process the closers and
transients one by one from left to right. Each time we connect the
vertex $i$ that we are processing to the $\xi_i$-th available
opener or transient to the left of $i$, where the openers and
transients are ranked from right to left. If we rank the openers
and transients from left to right, we get $\Phi_l((M,\xi))$. It
can be readily checked that $\Phi_r$ and $\Phi_l$ are well defined. 
Moreover:

\begin{proposition}
The maps $\Phi_r, \Phi_l:\Gamma_n \rightarrow \Pi_n$ are
bijections.
\end{proposition}
The proof can be found in \cite{flajolet, kz} and their references.

\begin{example}
If $(M,\xi)$ is the Charlier diagram in Figure~\ref{fig:charlier},
then
\begin{align*}
\Phi_r((M,\xi))&=\{\{1,7,10\},\{2,4,6,8\},\{3\},\{5,9\},\{11,12\}\}
\\
\Phi_l((M,\xi))&=\{\{1,4,6,7,9\},\{2,10\},\{3\},\{5,8\},\{11,12\}\}
\end{align*}
\begin{figure}[ht]
\begin{center}
\includegraphics[width=11cm]{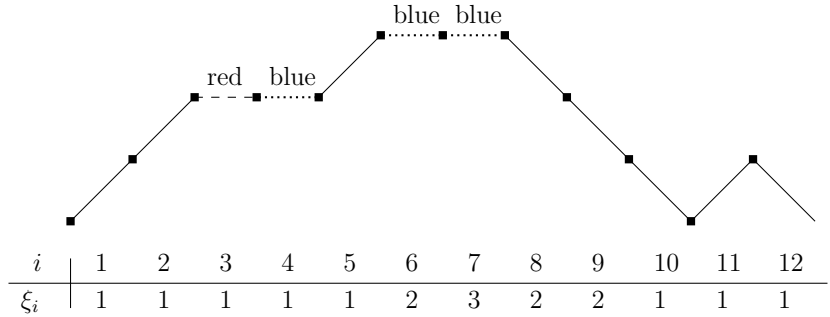}
\end{center}
\caption{ A Charlier diagram}
\label{fig:charlier}
\end{figure}
\end{example}

For $M \in RBM_n$  let $d_i$ be the number of NE and RE steps that
start at height $i$, $(i \geq 0)$. The \emph{profile} of $M$ is the
sequence $pr(M)=(d_0, \dots, d_l)$, where $l=\max\{i: d_i \neq 0
\}$. Note that this implies that $d_i \geq 1$ for each $ i=0,  \dots, l$,
and that the path $M$ is of height $l$ or $l+1$. The \emph{semi-type} of
$M=(s_1, \dots, s_i)$ is the sequence $st(M)=(\epsilon_1, \dots,
\epsilon_n)$ where $\epsilon_i=0$ if $s_i$ is a NE or RE step, and
$\epsilon_i=1$ if $s_i$ is a SE or BE step. For example, if  $M$
is the path in Figure \ref{fig:charlier}, then 
$pr(M)=(2, 1, 2)$, and $st(M)=(0, 0, 0, 1, 0, 1, 1, 1, 1, 1, 0, 1)$.  

Let $\lambda \in
\Pi_n$ and $\Phi^{-1}_{r}(\lambda)=(M,\xi^{r})$,
$\Phi^{-1}_{l}(\lambda)=(M,\xi^{l})$. 
Define $\varphi(\lambda)=M \in RBM_n$.
Note that for a given $\lambda$, $\varphi(\lambda)$ can be easily
constructed using the four steps above.
The next lemma gives the relation between a partition and its
corresponding restricted bicolored Motzkin path and Charlier
diagram.
%%%%%%%%%%%%%%%%%%%%%%%%%%%%%%%%%%%%%%%%%%%%%%%%%%%%%%%%%%%%%%%%%%%%%%%%%%%%%%%
% Lemma: properties of \Phi and \varphi                                       %
%%%%%%%%%%%%%%%%%%%%%%%%%%%%%%%%%%%%%%%%%%%%%%%%%%%%%%%%%%%%%%%%%%%%%%%%%%%%%%%
\begin{lemma}\label{L:properties}
Let $\Phi_r$, $\Phi_l$ and $\varphi$ be the maps defined
above and $\Phi^{-1}_{r}(\lambda)=(M,\xi^{r})$,
$\Phi^{-1}_{l}(\lambda)=(M,\xi^{l})$.
\begin{itemize}
 \item[(a)] The number of blocks of $\lambda$ is equal to the
total number of NE and RE steps of $M$.

\item[(b)] $cr(\lambda)=\sum_{i=1}^{n}{(\xi_i^r - 1)}$,
$ne(\lambda)=\sum_{i=1}^{n}{(\xi_i^l - 1)}.$

\item[(c)] $pr(M)=(d_0,\dots,d_l)$ \;if and only if \;
$crset(\lambda)=\{0^{d_0},\dots,l^{d_l}\}.$

\item[(d)] $neseq(\lambda)= v_1 \dots v_k$  if and only if the
zeros in $st(M)=(\epsilon_1, \dots, \epsilon_n)$ are in the
positions $v_1+1, v_2+2, \dots, v_k +k$. 
\end{itemize}
\end{lemma}

\begin{proof}
(a) The result follows from the fact that the number of blocks of
$\lambda$ is equal to the total number of openers and singletons.

(b) Denote by $E$ be the set of arcs of $\lambda$. For $e=(i,j)
\in E$ let $c_e = |\{(p,q) \in E : i < p < j <q\}|$. Then
$cr(\lambda)=\sum_{e \in E}{c_e }$. Similarly, if $n_e = |\{(p,q)
\in E : p < i < j <q\}|$, then $ne(\lambda)=\sum_{e \in E}{n_e }$.
{From} the definitions of $\Phi_r$ and $\Phi_l$ it follows that
$c_{(i,j)}= \xi_j^r-1$ and $n_{(i,j)}= \xi_j^l-1$. Hence the
claim.

(c) Using the notation at the beginning of Section~\ref{S:main},
we have $crseq(\lambda)=(u_1,\dots,u_k)$, where $u_i$ is
the number of edges $(p,q)$ such that $p < b_i <q$. Here 
$b_i=\min{ B_i}$, that is, $b_i$ is the $i$-th opener or singleton from
left to right. But the step in $M$ which corresponds to $b_i$ is
of height $h$ if and only if $u_i=h$.

(d) It follows directly from the definitions of $neseq(\lambda)$ and
$st(M)$.
\end{proof}
%%%%%%%%%%%%%%%%%%%%%%%%%%%%%%%%%%%%%%%%%%%%%%%%%%%%%%%%%%%%%%%%%%%%%%%%%%%
A composition of $k$ is an ordered tuple $(d_0,\dots,d_l)$ of
positive integers whose sum is $k$.

%%%%%%%%%%%%%%%%%%%%%%%%%%%%%%%%%%%%%%%%%%%%%%%%%%%%%%%%%%%%%%%%%%%%%%%%%%%%%%%
% Lemma2 in this section                                                      %
%%%%%%%%%%%%%%%%%%%%%%%%%%%%%%%%%%%%%%%%%%%%%%%%%%%%%%%%%%%%%%%%%%%%%%%%%%%%%%%
\begin{lemma}\label{L:cross}
 Let $l \geq 0, k \geq 1$, and $n \geq k$.
\begin{itemize}
\item[(a)]  If $\lambda \in \Pi_{n,k}$ and
$crset(\lambda)=\{0^{d_0},\dots,l^{d_l}\}$, then $(d_0,\dots,d_l)$
is a composition of $k$ into $l +1$ parts, where $l \leq n-k$,
 and $0 \leq cr(\lambda) \leq
 (n-k-1)l- \frac{l(l-1)}{2} $.

\item[(b)] Given a composition $(d_0,\dots,d_l)$ of $k$ into
$l+1 \leq n-k+1$ parts and an
integer $c$ such that $0 \leq c \leq (n-k-1)l- \frac{l(l-1)}{2}$,
there exists $\lambda \in \Pi_{n,k}$ with
$crset(\lambda)=\{0^{d_0},\dots,l^{d_l}\}$ and $cr(\lambda)=c$.
\end{itemize}

\end{lemma}

\begin{proof}
(a) It is clear that $d_0 + \cdots + d_l =k$. It follows that all
the $d_i$'s are positive from part (c) of
Lemma~\ref{L:properties}. Moreover, $\lambda$ has at least $l$
openers and, therefore, at least $l$ closers. So, $k+l \leq n$,
i.e., $l+1 \leq n-k+1$. Let $c_i$ (respectively $t_i$) be the
number of SE (respectively BE) steps at level $i$, $1 \leq i \leq
l+1$. Then $\sum_{i=1}^{l+1}{(c_i+t_i)}=n-k$ and $c_i \geq 1$, $1
\leq i \leq l$. Using part (b) of Lemma~\ref{L:properties}, we
have $0 \leq cr(\lambda) \leq
\sum_{i=1}^{l}{(1+0)(i-1)}+(n-k-l)l=(n-k-1)l- \frac{l(l-1)}{2}$.

(b) Suppose first that $l+1 \leq n-k$. Let $M \in RBM_n$ consist
of $d_0-1$ RE steps followed by a NE step, then $d_1-1$ RE steps
followed by one NE step, etc., $d_{l}-1$ RE steps followed by a NE
step, then $n-k-l-1$ BE steps, and $l+1$ SE steps. It is not hard
to see that indeed $M \in RBM_n$. The path never crosses the
$x$-axis and all the BE steps, if any, are at hight $l+1 \geq 1$.
Also $pr(M)=(d_0,\dots,d_l)$. Consider all the sequences
$\xi=(\xi_1,\dots,\xi_n)$ such that $(M,\xi)$ is a Charlier
diagram. Then

\begin{align}\label{E:xi}
\xi_i=1, \qquad & 1 \leq i \leq k \notag \\
1 \leq \xi_i \leq l+1, \qquad & k+1 \leq i \leq n-l-1 \\
1 \leq \xi_{n-i+1} \leq i, \qquad & 1 \leq i \leq l+1. \notag
\end{align}
Hence
\begin{eqnarray} \label{E:ineq1}
0 \leq \sum_{i=1}^{n}{(\xi_i-1)} \leq (n-k-l-1)l+l+(l-1)+
\cdots + 1= (n-k-1)l- \frac{l(l-1)}{2}.
\end{eqnarray}
In the case $l=n-k$, construct $M \in RBM_n$ similarly: $d_0-1$ RE
steps followed by a NE step, then $d_1-1$ RE steps followed by one
NE step, etc., $d_{l}$ RE steps, followed by  $l$ SE steps. (Note
that, unlike in the case $l < n-k$, the path $M$ is of height
$l$) All the sequences $\xi=(\xi_1,\dots,\xi_n)$ such that
$(M,\xi)$ is a Charlier diagram satisfy the following properties:
\begin{align}\label{E:xi1}
\xi_i=1, \qquad & 1 \leq i \leq k \notag \\
1 \leq \xi_{n-i+1} \leq i, \qquad & 1 \leq i \leq l. 
\end{align}
Hence \begin{eqnarray} \label{E:ineq2}
0 \leq \sum_{i=1}^{n}{(\xi_i-1)} \leq (l-1)+
\cdots + 1= (n-k-1)l- \frac{l(l-1)}{2}.
\end{eqnarray}
Because of~\eqref{E:ineq1} (respectively~\eqref{E:ineq2}), for any
integer $c$ between $0$ and $(n-k-1)l- \frac{l(l-1)}{2}$,
$\xi$ can be chosen to satisfy the conditions~\eqref{E:xi} (respectively~\eqref{E:xi1}) and such that
$\sum_{i=1}^{n}{(\xi_i-1)}=c$. Since $\Phi$ is a bijection, there is
$\lambda \in \Pi_{n,k}$ such that $\Phi(\lambda)=(M,\xi)$ and, by
part (b) and (c) of Lemma~\ref{L:properties}, $cr(\lambda)=c$
and $crset(\lambda)=\{0^{d_0}, \dots, l^{d_l}\}$.
\end{proof}
%%%%%%%%%%%%%%%%%%%%%%%%%%%%%%%%%%%%%%%%%%%%%%%%%%%%%%%%%%%%%%%%%%%%%%%%%%%%%%%%
%%%%%%%%%%%%%%%%%%%%%%%%%%%%%%%%%%%%%%%%%%%%%%%%%%%%%%%%%%%%%%%%%%%%%%%%%%%%%%%%
% Theorem: Number of crossing classes                                          %
%%%%%%%%%%%%%%%%%%%%%%%%%%%%%%%%%%%%%%%%%%%%%%%%%%%%%%%%%%%%%%%%%%%%%%%%%%%%%%%%
%%%%%%%%%%%%%%%%%%%%%%%%%%%%%%%%%%%%%%%%%%%%%%%%%%%%%%%%%%%%%%%%%%%%%%%%%%%%%%%%
\begin{theorem} Let $n \geq k \geq 1$ and $m=\min{\{n-k,k-1\}}.$ Then
\begin{eqnarray}\label{E:eq1}
|\Pi_{n,k}/\sim_{cr}|= \sum_{l=0}^{m}{\binom{k-1}{l}[(n-k-1)l-\frac{l(l-1)}{2}+1]}.
\end{eqnarray}
In particular, if $n \geq 2k-1$,
\begin{eqnarray}\label{E:eq2}
\left|\Pi_{n,k}/\sim_{cr} \right|=
(n-k-1)(k-1)2^{k-2}+2^{k-1}-(k-1)(k-2)2^{k-4}. 
\end{eqnarray}
\end{theorem}

\begin{proof}
Recall that $\lambda \sim_{cr} \pi$ if and only if $cr(\lambda)=cr(\pi)$ and
$crset(\lambda)=crset(\pi)$. Therefore,
$|\Pi_{n,k}/\sim_{cr}|=|\{(crset(\lambda),cr(\lambda)): \lambda \in \Pi_{n,k}\}|$.
 Using Lemma~\ref{L:cross} and the fact that the number of compositions of
$k$ into $l+1$ parts, $0 \leq l \leq k-1$, is $\binom{k-1}{l}$, we derive~\eqref{E:eq1}.
In particular, when $n \geq 2k-1$,
\[|\Pi_{n,k}/\sim_{cr}|= 
\sum_{l=0}^{k-1}{\binom{k-1}{l}[(n-k-1)l-\frac{l(l-1)}{2}+1]}.\]
 But
 \begin{align*}
 \sum_{l=0}^{k-1}{ \binom{k-1}{l}}&=(1+x)^{k-1}|_{x=1}=2^{k-1},\\
\sum_{l=0}^{k-1}{l \binom{k-1}{l}}&=
\left(\frac{d}{dx}(1+x)^{k-1}\right)\vert_{x=1}=(k-1)(1+x)^{k-2}|_{x=1}
\\&=(k-1)2^{k-2},
\\ \sum_{l=0}^{k-1}{l(l-1)
\binom{k-1}{l}}&=\left(\frac{d^2}{dx^2}(1+x)^{k-1}\right)\vert_{x=1}
\\ &=(k-1)(k-2)(1+x)^{k-3}|_{x=1}=(k-1)(k-2)2^{k-3},
\end{align*}
and~\eqref{E:eq2} follows.
\end{proof}

%%%%%%%%%%%%%%%%%%%%%%%%%%%%%%%%%%%%%%%%%%%%%%%%%%%%%%%%%%%%%%%%%%%%%%%%%%%%%%%%%%%%%%%%%%
Theorem 4.5 implies that there are many more examples of
different partitions $\lambda$ and $\pi$ for  which the 
statistic $cr$ has same distribution on the the levels of
$\mathcal{T}(\lambda)$ and $\mathcal{T}(\pi)$. For example,
$|\Pi_{2k,k}| > (2k-1)!! \approx \sqrt{2} \left( \frac{2k}{e}
\right)^{k}$ while $\left| \Pi_{2k,k}/\sim_{cr} \right|\approx
3k^{2}2^{k-4}$.

Next we analyze the number of nesting-similarity classes.  
%In~\cite{klazar}, Klazar uses Dyck
%paths to count the nesting-similarity classes. For partitions we
%need to work with restricted bicolored Motzkin paths. The key
%difference, which makes our computations more complicated, is the
%following: while $pr(D)$ determines the number of SE steps of the
%Dyck path $D$ at each level,  there exist Motzkin paths with same
%profile, but different total number of SE and BE steps at each
%level. ?? Hence the nesting-equivalence classes can not be counted
%as a summation over (restricted) Motzkin paths. This presented
%difficulty when trying to compute the number of nesting-similarity
%classes for set partitions. 
First we derive a recurrence for the
numbers $f_{n,k}=|\Pi_{n,k}/\sim_{ne}|$.

%%%%%%%%%%%%%%%%%%%%%%%%%%%%%%%%%%%%%%%%%%%%%%%%%%%%%%%%%%%%%%%%%%%%%%%%%%%%%%%%%%%%%%%%%%%%%%%
%  Theorem: Number of Nesting classes                                                         %
%%%%%%%%%%%%%%%%%%%%%%%%%%%%%%%%%%%%%%%%%%%%%%%%%%%%%%%%%%%%%%%%%%%%%%%%%%%%%%%%%%%%%%%%%%%%%%%
\begin{theorem} Let $n \geq k \geq 1$. Then
\begin{align}
f_{n,1}&=1, \label{E:nesting1}\\
f_{n,k}&=\sum_{r=k-1}^{n-1}{f_{r,k-1}} + (k-1)\binom{n-2}{k},
\qquad k \geq 2. \label{E:nesting2}
\end{align}
\end{theorem}

\begin{proof}
Equation \eqref{E:nesting1} is clear since $|\Pi_{n,1}|=1$.

Recall that $\lambda \sim_{ne} \pi$ if and only if
$ne(\lambda)=ne(\pi)$ and $neseq(\lambda)=neseq(\pi)$. By
Lemma~\ref{L:properties}, $f_{n,k}$ is equal to the number of
pairs $(\epsilon,c)$ such that there exists $\lambda \in
\Pi_{n,k}$ with $ne(\lambda)=c$ and
$st(\varphi(\lambda))=\epsilon$. It is not hard to see that for a
given a sequence $\epsilon =(\epsilon_1, \dots, \epsilon_n) \in
\{0,1\}^n$, there exists $\lambda \in \Pi_{n,k}$ such that
$st(\varphi(\lambda))=\epsilon$ if and only if $\epsilon$ has $k$
zeros and $\epsilon_1=0$. Denote the set of all such sequences by
$S^0_{n,k}$ and denote the set of all $\epsilon
 \in \{0,1\}^n$ with $k$ zeros by $S_{n,k}$.

For a sequence $\epsilon  \in S_{n,k}$  define a bicolored Motzkin
path $M=M(\epsilon)=(s_1, \dots, s_n)$ as follows. For $i$ from
$n$ to $1$ do:
\begin{enumerate}
\item[-] If $\epsilon_i=0$ and $s_i$ is not defined yet, then set
$s_i$ to be a RE step;

\item[-] If $\epsilon_i=1$ and there is $j<i$ such that
$\epsilon_j=0$ and $s_j$ is not defined yet, then set $s_i$ to be
a SE step and $s_{j_0}$ to be a NE step, where $j_0=\min\{j :
\epsilon_j=0 \text{\; and $s_j$ is not defined yet}\}$;

\item[-] If $\epsilon_i=1$ and there is no $j<i$ such that
$\epsilon_j=0$ and $s_j$ has not been defined yet, set $s_i$ to be
a BE step. 
\end{enumerate}
Note that we build $M$ backwards, from $(n,0)$ to $(0,0)$.
Let $h_i$ be the height of $s_i$ and
$ne(\epsilon)=\sum(h_i-1)$, where the sum is over all the indices
$i$ such that $\epsilon_i=1$.
For example, if $\epsilon=(0,0,0,1,0,1,1,1,1,1,0,1)$, then
\[
M(\epsilon)=(NE,NE,NE,BE,NE,BE,BE,SE,SE,SE,RE,SE).
\] 
The sequence of the heights of all the steps of $M$ is 
$(0,1,2,3,3,4,4,4,3,2,1,1)$ and
$ne(\epsilon)=(3-1)+(4-1)+(4-1)+(4-1)+(3-1)+(2-1)+(1-1)=14$.

Although clearly $M(\epsilon)$ stays above the $x$-axis, it is not
necessarily a restricted bicolored Motzkin path. The reason is that
any  $1$ in $\epsilon$ before the first zero would  produce a BE
step on the $x$-axis. Hence, $M(\epsilon) \in RBM_n$ if and only if
$\epsilon_1=0$, or equivalently,  $\epsilon \in S^0_{n,k}$.

We claim that for a fixed $\epsilon \in S^0_{n,k}$, there is
$\lambda \in \Pi_{n,k}$ such that $st(\varphi(\lambda))=\epsilon$
and $ne(\lambda)=c$ if and only if $0 \leq c \leq ne(\epsilon)$.
To show the if part, one can choose a sequence $\xi=(\xi_1, \dots,
\xi_n)$ with $1 \leq \xi_i \leq h_i$ if $\epsilon_i=1$, 
$\xi_i=1$ if $\epsilon_i=0$, 
and $\sum_{i=1}^{n}
(\xi_i-1)=c$. Then Lemma~\ref{L:properties} implies that
$\Phi_l((M,\xi))$ satisfies the requirements. Conversely, suppose
$\lambda \in \Pi_{n,k}$ is such that
$st(\varphi(\lambda))=\epsilon$. Let
$\Phi_{l}^{-1}(\lambda)=(M',\xi')$. Then the height $h_{i}^{'}$ of
each BE and SE step of $M'$ satisfies
\[ h_{i}^{'} \leq \min\{\text{\# zeros in \;} (\epsilon_1,\dots,\epsilon_{i-1}), (\text{\#
ones in \;} (\epsilon_{i+1},\dots,\epsilon_n)) +1\}=h_i. \] 
Now, by Lemma~\ref{L:properties}, 
\[
ne(\lambda) = \sum(\xi_{i}^{'}-1) \leq
\sum(h_{i}^{'}-1) \leq \sum(h_{i}-1) = ne(\epsilon).\] 
The claim is proved. Back to the proof of Theorem 4.6, we have
\[
f_{n,k}=\sum_{\epsilon \in S^0_{n,k}}(ne(\epsilon)+1)=\sum_{\epsilon \in S^0_{n,k}}{(\sum
(h_i-1)+1)}=\sum_{\epsilon \in
S^0_{n,k}}{\sum{h_i}}-(n-k-1)\binom{n-1}{k-1}.
\] 
Set
\[  g_{n,k}=\sum_{\epsilon \in
S^0_{n,k}}{\sum{h_i}} \qquad \text{and} \qquad
g^{*}_{n,k}=\sum_{\epsilon \in S_{n,k}}{\sum{h_i}},\] 
where the inner
sums are taken over all the indices $i$ such that $\epsilon_i=1$.
With this notation, 
\begin{equation}\label{E:fgrelation}
f_{n,k}=g_{n,k}-(n-k-1)\binom{n-1}{k-1}. 
\end{equation}
The sequences 
$g_{n,k}$ and $g^{*}_{n,k}$ satisfy the following recurrence
relations:
\begin{align}
g_{n,k}&=g_{n-1,k-1}+g^{*}_{n-2,k-1}+(n-k)\binom{n-2}{k-1}, \label{E:g1}\\
 g_{n,k}^{*}&=\sum_{r=k}^{n}g_{r,k}. \label{E:g2}
\end{align}
To see~\eqref{E:g1}, note that if $\epsilon_n=0$ then
$(\epsilon_1,...,\epsilon_{n-1}) \in S^0_{n-1,k-1}$ and
$M(\epsilon)$ is $M(\epsilon_1,...,\epsilon_{n-1})$ with one RE
step appended, and if $\epsilon_n=1$ then
$(\epsilon_2,...,\epsilon_{n-1}) \in S_{n-2,k-1}$ and
$M(\epsilon_2,...,\epsilon_{n-1})$ is obtained from $M(\epsilon)$
by deleting the first NE and the last SE step. For~\eqref{E:g2},
if $\epsilon_1=\cdots =\epsilon_{r-1}=1$ and $\epsilon_{r}=0$,
then $M(\epsilon_r, \dots , \epsilon_n)$ is obtained from
$M(\epsilon)$ by deleting the first $r-1$ BE steps at level 0.
Substituting~\eqref{E:g2} into~\eqref{E:g1} gives
\begin{equation}\label{E:grec}
g_{n,k}=\sum_{r=k-1}^{n-1}g_{r,k-1} + (n-k)\binom{n-2}{k-1}.
\end{equation}
Finally, by substituting $g_{n,k}$ from~\eqref{E:fgrelation}
into~\eqref{E:grec} and simplifying, we obtain~\eqref{E:nesting2}.
\end{proof}

\begin{corollary}

\begin{align*}
|\Pi_1/\sim_{ne}|&=1, \qquad |\Pi_2/\sim_{ne}|=2 \\
|\Pi_n/\sim_{ne}|&=2^{n-5}(n^2-5n+22), \qquad n \geq 3
\end{align*}
\end{corollary}

\begin{proof}
Denote $|\Pi_n/\sim_{ne}|$ by $F_n$. Using
$F_n=\sum_{k=1}^{n}f_{n,k}$,~\eqref{E:nesting1},
and~\eqref{E:nesting2}, we get
\begin{align*} F_n&=1+F_1+ \cdots + F_{n-1} + \sum_{k=2}^{n}(k-1)\binom{n-2}{k} \\
                   &=F_1+ \cdots +F_{n-1}+ (n-4)2^{n-3}+2, 
\qquad \qquad n \geq 2.
\end{align*}
This yields the recurrence relation
\[ F_n=2F_{n-1} + (n-3)2^{n-4}, \qquad  n\geq 3 \]
with initial values $F_1=1$ and $F_2=2$, 
which has the solution
%\[ F_n=1+F_1+ \cdots + F_{n-1} + \sum_{k=2}^{n}(k-1)\binom{n-2}{k}
%                   =F_1+ \cdots +F_{n-1}+ (n-4)2^{n-3}+2 \]
\[F_n=2^{n-5}(n^2-5n+22), \qquad  n \geq 3.\]
\end{proof}

The following tables give the number of crossing/nesting-similarity
classes on $\Pi_{n,k}$ for small $n$ and $k$. 

\begin{center}
\parbox{6cm}{
\begin{tabular}{c|cccccc}
\multicolumn{7}{c}{crossing-similarity classes} \\ \hline
$n \setminus k $ & 1 & 2  & 3 &4  & 5 & 6  \\ \hline
1                & 1 &    &   &   &   &    \\
2                & 1 & 1  &   &   &   &    \\
3                & 1 & 2  & 1 &   &   &    \\
4                & 1 & 3  & 3 & 1 &   &    \\
5                & 1 & 4  & 7 & 4 & 1 &    \\
6                & 1 & 5 & 11 & 4 & 5 & 1  \\ 
\end{tabular}
}
\parbox{6cm}{
\begin{tabular}{c|cccccc}
\multicolumn{7}{c}{nesting-similarity classes} \\ \hline
$n \setminus k $ & 1 & 2  & 3 &4  & 5 & 6  \\ \hline
1                & 1 &    &   &   &   &    \\
2                & 1 & 1  &   &   &   &    \\
3                & 1 & 2  & 1 &   &   &    \\
4                & 1 & 4  & 3 & 1 &   &    \\
5                & 1 & 7  & 9 & 4 & 1 &    \\
6                & 1 & 11 & 22 & 16 & 5 & 1 
\end{tabular}
} 
\end{center}

The two equivalence relations $\sim_{cr}$ and $\sim_{ne}$ on set partitions 
are not compatible. 
From the tables it is clear that $\sim_{cr}$ is not a refinement 
of $\sim_{ne}$. On the other hand, 
let $\pi=\{ \{1,3\}, \{2,4\}, \{5,6\} \}$ and
$\lambda=\{ \{ 1,3,6\}, \{2,4\}, \{5\} \}$. 
It is easy to check that  
 $\pi \sim_{ne} \lambda$,  
but  $\pi \not\sim_{cr} \lambda $,  as $cr(\pi)=1$ and 
$cr(\lambda)=2.$  

%%%%%%%%%%%%%%%%%%%%%%%%%%%%%%%%%%%%%%%%%%%%%%%%%%%%%%%%%%%%%%%%%%%%%%%%%%%%%%%%%%%%%%%%%%%%%%%%%%%%%
%%%%%%%%%%%%%%%%%%%%%%%%%%%%%%%%%%%%%%%%%%%%%%%%%%%%%%%%%%%%%%%%%%%%%%%%%%%%%%%%%%%%%%%%%%%%%%%%%%%%%
% Section: Generating function                                                                      %
%%%%%%%%%%%%%%%%%%%%%%%%%%%%%%%%%%%%%%%%%%%%%%%%%%%%%%%%%%%%%%%%%%%%%%%%%%%%%%%%%%%%%%%%%%%%%%%%%%%%%
\section{Generating function for crossings and nestings}

In this section we analyze the generating function \[
S_\pi(q,p,z)= \sum_{l \geq 0}\sum_{\lambda \in \mathcal{T} (\pi ,
l) }q^{cr(\lambda)} p^{ne(\lambda)}z^l \] for a given partition
$\pi$, and  derive a continued fraction expansion for
$S_\pi(q,p,z)$.
{For} this we work with the group $G=\mathbb{Z} \oplus \mathbb{Z}$
and $\alpha =(1,0), \beta=(0,1)$. Fix a partition $\pi$ with $k$
blocks. Define
$E_{\alpha,\beta}(\pi,l)=\cup_{m=k}^{k+l}{E_{\alpha,\beta}(\pi,l,m)}$,
i.e., $E_{\alpha,\beta}(\pi,l)$ is the multiset of sequences
$seq_{\alpha,\beta}{(\mu)}$ associated to the partitions $\mu \in
\mathcal{T}(\pi,l)$. A recurrence analogous to~\eqref{E:levelseq}
holds. Namely, for $l \geq 1$
\begin{equation}\label{E:levelseq2}
E_{\alpha,\beta}(\lambda,l)=
R_{\alpha,\beta}(E_{\alpha,\beta}(\lambda,l-1)) \cup
M(E_{\alpha,\beta}(\lambda,l-1)).
\end{equation}
For simplicity we write $E_l$
instead of $E_{\alpha,\beta}(\pi,l)$ when there is no confusion.
 Define $b_{l,r}$ to be the generating function of the
multiset $f_{0}^{r}(E_l)$, i.e., 
\[
b_{l,r}(q,p)=\sum_{(x,y) \in f_{0}^{r}(E_l)}q^x p^y,
\] 
where
$(x,y) \in f_{0}^{r}(E_l)$ contributes to the sum above according
to its multiplicity in $f_{0}^{r}(E_l)$. By convention, let
$b_{l,r}(q,p)=0$ if $f_{0}^{r}(E_l)=\emptyset$, or, one of
$l, r$ is negative.  For simplicity
we  write $b_{l,r}$ for $b_{l,r}(q,p)$. Note that
$b_{l,0}=\sum_{\lambda \in \mathcal{T} (\pi , l) }q^{cr(\lambda)}
p^{ne(\lambda)}$ and hence 
\begin{eqnarray} \label{S_and_b} 
S_\pi(q,p,z)= \sum_{l \geq 0}b_{l,0}z^l.
\end{eqnarray} 
By the formulas ~\eqref{E:levelseq2}, ~\eqref{E:mrec}, and
the proof of part (c) of Lemma~\ref{L:big},  we get
\begin{align*} f_0^r(E_l) &=f_0^r(M(E_{l-1})) \cup f_0^r(R_{\alpha,
\beta}(E_{l-1}))\\
&=f_0^{r-1}(E_{l-1}) \cup f_0^r(E_{l-1}) \cup
\bigcup_{t=0}^{r}f^{r+1}_{t\alpha+(r-t)\beta}(E_{l-1})\cup
\bigcup_{t=0}^{r}f^{r}_{t\alpha+(r-t)\beta}(E_{l-1}), 
\end{align*}
which leads to a recurrent relation for $b_{l,r}$:
\[b_{l,r}= b_{l-1,r-1}+b_{l-1,r}+(\sum_{t=0}^{r}q^t
p^{r-t})b_{l-1,r+1}+(\sum_{t=0}^{r}q^t p^{r-t})b_{l-1,r}.\] Using
the standard notation $[r]_{q,p}:=\frac{q^r-p^r}{q-p}$, we can
write this as
\begin{proposition}\label{L:rec}
\[b_{l,r}=
b_{l-1,r-1}+(1+[r+1]_{q,p})b_{l-1,r}+[r+1]_{q,p}b_{l-1,r+1}. \]
\end{proposition}
If the sequence associated to the partition $\pi$ is $x_1x_2 \dots
x_k$, with $x_i=u_i \alpha +v_i \beta, \;1 \leq i \leq k$, then
\begin{eqnarray} \label{b_{0,r}}
b_{0,0}&=&q^{u_{1}}p^{v_{1}} \notag \\
b_{0,r}&=& \sum_{1< i_1 < \cdots <i_r \leq
k}q^{u_{i_{1}}+\cdots +u_{i_{r}}-(r-1)u_1}p^{v_{i_{1}}+\cdots
+v_{i_{r}}-(r-1)v_1} \quad \text{for} \; r \geq 1.
 \end{eqnarray}
In particular, $b_{0,r}=0$ if $r \geq k$.

Given $l$ and $s$, nonnegative integers, consider the paths from
$(l,0)$ to $(0,s)$ using steps $(-1,0)$, $(-1,1)$, and $(-1,-1)$
which do not go below the $x$-axis. Each step $(-1,0)$ \:
($(-1,1)$,\; $(-1,-1)$  respectively) starting at the
line $y=r$ has weight \: $[r+1]_{q,p}$ \; ($1+[r+1]_{q,p}$, \; $1,
\; \text{respectively})$. The weight $w(M)$ of such a path $M$ is
defined to be the product of the weights of its steps. Let
$c_{l,s}=\sum{w(M)}$, where the sum is over all the paths $M$
described above. Then from  Proposition~\ref{L:rec} one has
\[b_{l,0}=\sum_{0 \leq s \leq k-1} c_{l,s}b_{0,s}. \]
Set $a_{r}=[r+1]_{q,p}$ and $c_{r}=[r+1]_{q,p}+1$. By the
well-known theory of continued fractions (see~\cite{flajolet}),
$c_{l,s}$ is equal to the coefficient in front of $z^l$ in
\begin{eqnarray} \label{K_s}
K_s(z):=J^{/0/}(z)a_{0}zJ^{/1/}(z)a_{1}z \cdots
J^{/s/}(z)=\frac{1}{z^s}(Q_{s-1}(z)J(z)-P_{s-1}(z))
\end{eqnarray}  
where
\[J^{/h/}(z)=\cfrac{1}{1-c_{h}z-\cfrac{a_{h}z^2}{1-c_{h+1}z-\cfrac{a_{h+1}z^{2}}{\ddots}}}\]
and $\frac{P_k(z)}{Q_{k}(z)}$ is the $k$-th convergent of
$J(z):=J^{/0/}(z)$. 
Hence 
\begin{theorem} \label{gf:cf}
Let $\pi$ be a partition with $k$ blocks whose associated sequence is 
$x_1x_2 \dots x_k$, where $x_i=u_i \alpha + v_i \beta$ for $1 \leq i
\leq k$. Then 
\[
S_\pi(q, p,z) = \sum_{0 \leq s \leq k-1} b_{0,s} K_s(z), 
\]
where $b_{0,s}$ is given by the formula \eqref{b_{0,r}}, and $K_s(z)$
is given  by \eqref{K_s}.
\end{theorem} 

In particular, when $k=1$, i.e., $\pi$ is a
partition with only one block, then $b_{0,0}=1$ and 
$b_{l,0}=c_{l,0}b_{0,0}=c_{l,0}$.
Therefore
\begin{corollary} \label{T:sec}
If $|\pi|=1$, then
\[S_{\pi}(q,p,z)=\cfrac{1}{1-([1]_{q,p}+1)z-\cfrac{[1]_{q,p}z^2}{1-([2]_{q,p}+1)z-\cfrac{[2]_{q,p}z^{2}}{\ddots}}}.\]
\end{corollary}

\textsc{Remark.} 
Corollary~\ref{T:sec} leads to  a continued fraction
expansion for the generating function of crossings and nestings over
$\Pi$: Just  taking $\pi$ to be the
partition of $\{1\}$, and bearing in mind
that we are counting the empty partition as well, we get
\begin{align}\label{E:fraction2}
\sum_{n \geq 0}\sum_{\lambda \in \Pi_n  }q^{cr(\lambda)}
p^{ne(\lambda)}z^n&= 1+zS_{\{1\}}(q,p,z)\notag \\&=
1+\cfrac{z}{1-([1]_{q,p}+1)z-\cfrac{[1]_{q,p}z^2}{1-([2]_{q,p}+1)z-\cfrac{[2]_{q,p}z^{2}}{\ddots}}}.
\end{align}
A different expansion was given in~\cite{kz}, as
\begin{eqnarray}\label{E:fraction1}
\sum_{n \geq 0}\sum_{\lambda \in \Pi_n  }q^{cr(\lambda)}
p^{ne(\lambda)}z^n=
\cfrac{1}{1-z-\cfrac{z^2}{1-([1]_{q,p}+1)z-\cfrac{[2]_{q,p}z^2}{1-([2]_{q,p}+1)z-\cfrac{[3]_{q,p}z^2}{\ddots}}}}.
\end{eqnarray}

The fractions~\eqref{E:fraction1} and~\eqref{E:fraction2} can be
transformed  into each another by applying  twice the following contraction
formula for continued fraction, (for example, see \cite{clarke}):
\[
\cfrac{c_0}{1-\cfrac{c_1z}{1-\cfrac{c_2z}{\ddots}}}=
c_0+\cfrac{c_0c_1z}{1-(c_1+c_2)z-\cfrac{c_2c_3z^2}{1-(c_3+c_4)z-\cfrac{c_4c_5z^{2}}{\ddots}}}.
\]

\end{document}